\documentclass{amsart}
\usepackage{amsmath,amssymb}
\newtheorem{theorem}{Theorem}[section]

\newtheorem{proposition}[theorem]{Proposition}
\newtheorem{corollary}[theorem]{Corollary}

\newtheorem{remark}[theorem]{Remark}
\newtheorem{example}[theorem]{Example}

\begin{document}

\title{The $\lambda$~- dimension of commutative arithmetic rings}
\author{Fran\c{c}ois Couchot}

\begin{abstract} It is shown that every commutative arithmetic ring $R$ has
$\lambda$-dimension $\leq 3$. An example of a commutative Kaplansky ring
with $\lambda$-dimension 3 is given. Moreover, if $R$ satisfies one of
the following conditions, semi-local,  semi-prime, self
$fp$-injective, zero-Krull dimensional, CF or FSI 
then $\lambda$-$dim (R) \leq 2$. It is also shown that every zero-Krull
dimensional commutative arithmetic ring is a Kaplansky ring and
an adequate ring, that every B\'ezout ring with compact minimal prime
spectrum is Hermite and that each B\' ezout fractionnally self
$fp$-injective ring is a Kaplansky ring.
\end{abstract}
\maketitle

\section{Introduction, definitions and preliminaries}
\label{S:pre} 
All rings in this paper are commutative with unity and modules are
unitary. Following P. V\'amos \cite{Vam}, if ${\mathcal P}$ is a ring
property, we say that a ring $R$ is {\it locally} ${\mathcal P}$ if
$R_M$ has  ${\mathcal P}$ for every maximal ideal $M$, and
$R$ is {\it fractionnally} ${\mathcal P}$ if the
classical quotient ring $Q(R/A)$ of $R/A$ has ${\mathcal P}$ for
every proper ideal $A$ of $R$.

An $R$-module $E$ is said to be of {\it finite
$n$-presentation} if there exists an exact sequence:

\centerline{$F_n \rightarrow
F_{n-1} \rightarrow\cdots F_1 \rightarrow F_0 \rightarrow E \rightarrow 0$}

with the $F_i$'s free $R$-modules of finite rank. We write 
$\lambda_R(E) = \sup \{ n \mid$ there is a finite $n$-presentation of
$E\}$. If $E$ is not finitely generated we also put $\lambda_R(E)= - 1$.

The {\it $\lambda$-dimension} of a ring $R$ ($\lambda$-$dim
(R))$ is the least integer $n$ (or $\infty$ if none such
exists) such that $\lambda_R(E) \geq n$ implies $\lambda_R(E) =
\infty$. See \cite[chapter 8]{Vas}. Recall that $R$ is noetherian if and
only if $\lambda$-$dim (R) = 0$ and $R$ is coherent if and only if
$\lambda$-$dim (R) \leq 1$. 

This notion of $\lambda$-dimension of a
ring was formulated in \cite[chapter 8]{Vas} to study the rings of
polynomials or power series over a coherent ring.

In section~\ref{S:dim} of this paper it is proved that every
arithmetic ring has a $\lambda$-dimension $\leq 3$. We
give an example of a Kaplanky ring whose the $\lambda$-dimension is
exactly $3$. However, if an arithmetic ring satisfies an additional
property,(reduced, self $fp$-injective, semi-local, CF or
fractionnally self-injective), its $\lambda$-dimension is at most $2$.

In section~\ref{S:frac} we study fractionnally self $fp$-injective
rings. We prove that every reduced factor ring of a fractionnally self
$fp$-injective ring is semihereditary. It is shown that each
fractionnally self $fp$-injective ring which is B\'ezout is
Kaplansky. To state this last result, we give a positive answer to a
question of Henriksen by proving that any B\'ezout ring
with compact minimal prime spectrum is Hermite.

An $R$-module $E$ is said to be {\it uniserial} if the set of
its submodules is totally ordered by inclusion. A ring $R$ is a
{\it valuation ring} if $R$ is a uniserial module, and $R$ 
is {\it arithmetic} if $R$ is locally a valuation ring. A ring
is a {\it B\'ezout ring} if every finitely generated ideal is
principal. A ring $R$ is an {\it Hermite ring} if $R$ 
satisfies the following property~: for every $(a,b) \in R^2$, there
exist $d, a', b'$ in $R$ such that $a = da'$, $b = db'$ and
$Ra' + Rb' = R$. We say that $R$ is a {\it Kaplansky
ring} (or an {\it elementary divisor ring}) if for every matrix $A$,
with entries in $R$, there  exist
a diagonal matrix $D$ and invertible matrices $P$ and $Q$, with
entries in $R$, such that $PAQ = D$.  Then we have the following
implications~:

\centerline{Kaplansky ring $\Rightarrow$ Hermite ring
$\Rightarrow$ B\'ezout ring $\Rightarrow$ arithmetic ring~;}

but these implications are not reversible \cite{GiHe1} or \cite{Car}.

Recall that  $R$ is a Kaplansky ring if and only if every
finitely presented module is a finite direct sum of cyclic finitely
presented modules (\cite{Kap} and \cite{LLS}). We say that $R$ is an
{\it adequate ring} if $R$ is a B\'ezout ring satisfying the
following property~: for every $(a,b) \in R^2$, $a\not= 0$, there exist
$r$ and $s$ in $R$ such that $a = rs$, $Rr+Rb = R$, 
and if $s'$ is a nonunit that divides $s$, then $Rs'+Rb \not= R$.
An exact sequence $0 \rightarrow F \rightarrow E \rightarrow G \rightarrow 0$  is {\it pure}
if it remains exact when tensoring it with any $R$-module. In this case
we say that $F$ is a pure submodule of $E$. When $R$ is an
arithmetical ring then $F$ is a pure submodule of  $E$ if and only
if $rF = rE\cap F$ for every $r\in R$, \cite[Theorem 3]{War}.

The following proposition will be useful to provide us many examples in
the second part of this paper.

\begin{proposition}
\label{P:E}
 Let $I$ be an infinite set, $D$ a valuation domain and $N$ its
 maximal ideal. We consider 
\[S = \{f : I \rightarrow D \mid f\hbox{ constant except on a finite
 subset of }I\}\]. 
Then the following statements are true.
\begin{enumerate}
\item $S$  is a free $D$-module with basis  ${\mathcal B} = \{
\mathbf{1}, \mathbf{e}_i \mid i\in I\}$ where $\mathbf{1}(j) = 1$, and
$\mathbf{e}_i(j) =\delta_{ij}$, for every $j\in I$, and where
$\delta_{ij}$ is the Kronecker symbol.
\item $M_0 = N\mathbf{1} + \sum_{i\in I} D\mathbf{e}_i$ and $M_i =
D(\mathbf{1}-\mathbf{e}_i) + N\mathbf{e}_i$, for every $i\in I$, are
the maximal ideals of $S$. Moreover, $S_{M_{0}} \simeq D$  and
$S_{M_{i}} \simeq D$ for every $i\in I$.  The Jacobson radical $J(S) =
SN = N\mathbf{1} + \sum_{i\in I} N\mathbf{e}_i$.
\item $S$ is a Kaplansky ring and an adequate ring.
\end{enumerate}
\end{proposition}

\textbf{Proof.}
\begin{enumerate}
\item This assertion is obvious.
\item Let $M$ be a maximal ideal of $R$. If $\mathbf{e}_i\in M$, for every 
$i\in I$, then $M = M_0$.  The ideal $E$ of $S$ generated by
$\{\mathbf{e}_i \mid i\in I\}$ is a pure ideal of $S$, hence $S/E$ is a
flat $S$-module and $S/E \simeq D$. From this we deduce that 
$S_{M_{0}} \simeq D$.  If there exists $i\in I$ such that $\mathbf{e}_i
\notin M$, then $(\mathbf{1}-\mathbf{e}_i) \in M$ and we have $M =
M_i$. Moreover, $S/S(\mathbf{1}-\mathbf{e}_i)$ is a projective 
$S$-module and $S/S(\mathbf{1}-\mathbf{e}_i)\simeq D$. We  deduce
that $S_{M_{i}} \simeq D$.  Now, it is easy to get that 
$J(S) = N\mathbf{1} + \sum_{i\in I} N\mathbf{e}_i$.
\item By using the basis ${\mathcal B}$ of $S$ over $D$, it is easy
to prove that  $S$ is an Hermite ring and an adequate ring. From 
\cite[Theorem 8]{GiHe2} we deduce that $S$ is a Kaplansky ring.
\qed
\end{enumerate} 

\bigskip
An $R$-module $E$ is {\it fp-injective} if 
$\hbox{Ext}_R^1 (F, E) = 0$  for any finitely presented $R$-module $F,$
and $R$ is {\it self fp-injective} if $R$ is 
fp-injective as $R$-module. Recall that a valuation ring $R$ is
self fp-injective if and only if the set $Z(R)$ of its zero
divisors is its maximal ideal, \cite[Theorem 2.8]{Cou}.  We recall that a
module $E$ is fp-injective if and only if it is a pure submodule
of every overmodule.

We denote respectively $Spec(R)$, $MaxSpec(R)$ and $MinSpec(R),$ the
space of prime ideals, maximal ideals, and minimal prime ideals of
$R$, with the Zariski topology. If  $X = Spec(R), MaxSpec(R)$ or $MinSpec
(R),$ and $A$ a subset of $R$, then we denote  $V(A) = \{ P\in X \mid A
\subseteq P\}$ and $D(A) = \{ P\in X \mid A\not\subseteq P\}$.

Finally if $E$ is an $R$-module, $flat$-$dim (E)$
is the least integer $n$ such that\\ $\hbox{Tor}_{n+1}^R (F, E) = 0$ 
for every  $R$-module $F,$ and\\ $gl$-$w$-$dim (R) = \sup \{
flat$-$dim (E) \mid  E  \ \  R$-module$\}$.

\section{The $\lambda$-dimension}
\label{S:dim}
We begin with the more general result of this part.

\begin{theorem} 
\label{T:dimain}
Let $R$ be an arithmetic ring. Then the
following statements are true.
\begin{enumerate}
\item $\lambda$-$dim (R) \leq 3$.
\item If $R$ is a reduced ring then  $\lambda$-$dim (R) \leq 2$.
\end{enumerate}
\end{theorem}

\textbf{Proof.}
\begin{enumerate}
\item Let $E$ be a module such that $\lambda_R(E)\geq 3.$ We
consider the following finite 3-presentation of $E$:
 \[ F_3\overset{u_3}{\rightarrow} F_2 \overset{u_2}{\rightarrow} F_1
 \overset{u_1}{\rightarrow} F_0\rightarrow E \rightarrow 0.\]
 We choose bases ${\mathcal
B}_0$ and ${\mathcal B}_1$ of $F_0$ and $F_1$ respectively, and
let $A$ be the matrix associated with $u_1,$ with respect to our
given bases. Let $M$ be a maximal ideal of $R$.  By \cite[Theorem
1]{War} $E_M$is  a direct sum of cyclic finitely presented $R_M$-modules.
Therefore there exist a diagonal matrix $D$ and two invertible
matrices $P$ and $Q,$ with entries in $R_M$ such that $PAQ
= D.$ It is not difficult to find $t\in R \setminus M,$ such that 
$P$ and $Q$ are invertible matrices with entries in $R_t,$ $D$
a diagonal matrix with entries in  $R$ such that  $PAQ = D.$
It follows that there exist $a_1,\ldots, a_n\in R$ such that
$E_t \simeq \bigoplus_{k=1}^n(R_t/a_kR_t).$ Since
$\lambda_{R_{t}} (E_t) \geq 2,$ we deduce that $(0 :_{R_t} a_k )$ is
a finitely generated ideal of $R_t,$ and there exists $b_k\in R,$
such that $(0 :_{R_M} a_k ) = b_k R_M,$ for every $k, 1\leq
k\leq n.$ By multiplying $t$ with an element of $R\setminus M,$
we may assume that $(0 :_{R_t} a_k) = b_kR_t$ for every $k, 1\leq
k\leq n.$ Now, since $\lambda_{R_{t}} (E_t) \geq 3,$ by the same
way, we get that there exists $c_k\in R,$ such that $(0 :_{R_t} b_k) =
c_kR_t$ for every $k$, $1\leq k\leq n$.  Then, the equality $(0 :_{R_t}
a_k) = b_kR_t$ implies that $(0 :_{R_t} c_k ) = b_k R_t,$ for every
$k$, $1\leq k\leq n.$ Hence $\lambda_{R_{t}} (E_t) \geq 4.$

If we denote $U_M = D(t),$ then $(U_M)_{M\in MaxSpec(R)}$ is an open
overing of $MaxSpec(R),$ and since this space is quasi-compact, a finite
number  of these open subsets cover $MaxSpec(R).$ Thus, $MaxSpec(R)=
\cup_{j=1}^m U_j,$ where $U_j = D(t_j).$ Let\\ $K = ker(u_3).$
Now, for every $j$, $1\leq j \leq m$, $K_{t_{j}}$ is a finitely
generated $R_{t_{j}}$-module, hence there exists a finite subset $G_j$
of $K$ such that $K_{t_{j}} = \sum_{g\in G_{j}}R_{t_j}g.$ Then  $K$ is
generated by $\cup_{j=1}^m G_j$ and we get that $\lambda_R
(E) \geq 4.$
\item If $R$ is reduced, then $R_M$ is a valuation domain for
every maximal ideal $M$ of $R.$ Consequently $gl$-$w$-$dim (R)
\leq 1,$ and from \cite[Chapter 8]{Vas} we deduce that $\lambda$-$dim (R)
\leq 2.$ We can also deduce this result from our following
Corollary~\ref{C:wdone}.
\qed
\end{enumerate}

\bigskip
The example 1.3b of \cite{Vas} is a reduced arithmetic ring of
$\lambda$-dimension 2. Now, to complete the proof of our
Theorem~\ref{T:dimain}, an example of arithmetic ring with
$\lambda$-dimension 3 must be given.

\begin{example}
\label{E:dim3} 
\textnormal{Let $S$ be the ring defined in Proposition~\ref{P:E}.  We
  suppose that $D$ has a nonzero and nonmaximal prime ideal $J.$ Let $a\in
N\setminus J, b\in J, \ b\not= 0$ and $A = Dab\mathbf{1} +
\sum_{i\in I}J\mathbf{e}_i.$ We denote $R = S/A$ and $\overline{r} = r+A$ for every $r\in S$.  Then $R$ is a Kaplansky
ring and also an adequate ring since $A\subseteq J(S)$ by
  \cite[Proposition 4.4]{LLS}. Now it is easy to prove that $(0 :
  \overline{a\mathbf{1}}) =R\overline{b\mathbf{1}}$ and $(0 : \overline{b\mathbf{1}}) = R\overline{a\mathbf{1}} +
  \sum_{i\in I}R\overline{\mathbf{e}}_i.$ We deduce from this that 
$\lambda_R(R/R\overline{a\mathbf{1}}) = 2.$ Hence $\lambda$-$dim (R) =
3.$}
\end{example}

\begin{example} 
\label{E:redu}
\textnormal{Let $S$ be the ring defined in Proposition~\ref{P:E}, $A =
\sum_{i\in I}N\mathbf{e}_i$ and $R = S/A$. If $M$ is a
maximal ideal of $S$, we denote $\overline{M} = M/A$.  Then it is
easy to prove that $R_{\overline{M}_{0}} \simeq D$ and 
$R_{\overline{M}_{i}} \simeq D/N$ for every 
$i\in I$.  Consequently $R$ is a reduced ring, a Kaplansky
ring and an adequate ring. For every $a\in N$,
$a\not= 0$, $(0 : \overline{a\mathbf{1}}) = \sum_{i\in I}R \ \overline{\mathbf{
e}}_i$ is not finitely generated. Then $\lambda$-$dim (R) = 2$.  When
$D = \mathbb{Z}_2$, the ring of 2-adics numbers, and $I = \mathbb{N}$,
we obtain the example 1.3b of \cite{Vas}, if, in this example we
replace $\mathbb{Z}$ with $\mathbb{Z}_2$.}
\end{example}

\begin{theorem}
\label{T:dinj} 
Let $R$ be an arithmetic ring.

 If $R$ is
self fp-injective then $\lambda$-$dim (R) \leq 2$.
\end{theorem}

\textbf{Proof.} 
Let $E$ be a module with $\lambda_R (E) \geq 2$.  As in
the proof of Theorem~\ref{T:dimain},  for every 
$M\in MaxSpec(R)$, we can find $t\in R\setminus M$, $a_1,\ldots,a_n\in
R$, $b_1,\ldots,b_n\in R$ such that  $E_t \simeq
\bigoplus_{k=1}^n(R_t/a_kR_t)$
 and $(0 :_{R_t} a_k) = b_kR_t$ for  every $k$, $1\leq k\leq n$.

Since $\lambda_R(E) \geq 2$, then the canonical homomorphism\\ \(\bigl(
\hbox{Ext}_R^1 (E, R)\bigr)_t \rightarrow \hbox{Ext}_{R_{t}}^1 (E_t,
R_t)\) is an isomorphism.  Thus\\ $\hbox{Ext}_{R_{t}}^1 (E_t, R_t) = 0$ and 
$\hbox{Ext}_{R_{t}}^1 ( R_{t}/a_{k}R_{t},
R_t) = 0$\\  for every $k,$ $1\leq k \leq n$.  From the
following projective resolution of  $R_{t}/a_{k}R_{t}$ :
\(R_t\overset{b_k}{\rightarrow} R_t\overset{a_k}{\rightarrow} R_t\),
we deduce that $(0 :_{R_t}b_k) = a_k R_t$ \ for every $k$, $1\leq
k\leq n$.  Hence $\lambda_{R_{t}} (E_t) \geq 3$.  Now, as in the proof
of Theorem~\ref{T:dimain}, we get that $\lambda_R (E) \geq 3$.
\qed

\bigskip

When $R$ is a reduced ring we have a more general result.

\begin{theorem} 
\label{T:redu}
Let $R$ be  a reduced ring. Then $R$ is
self fp-injective if and only if $R$ is a Von Neumann regular
ring.
\end{theorem}

\textbf{Proof.}
Only necessity requires a proof. Since $R$ is reduced, $R$ is a subring of $S =
\Pi_{P\in MinSpec(R)}Q(R/P)$, and $S$ is a Von Neumann
regular ring. Hence, for every $r\in R$, there exists $s\in S$ 
such that $r^2s = r$.  But, since $R$ \ is self $fp$ - injective, 
$R$ is a pure submodule of $S$.  Thus, there exists $s'\in R$ such
that $r^2s' = r$.
\qed

\begin{remark}
\textnormal{We can prove that for every $n\in \mathbb{N}$, there
  exists a self injective ring $R$ such that $\lambda$-$dim (R) = n$.
  Let  $D$ be a local noetherian regular ring, $N$ its maximal ideal
  and $E$ the $D$-injective hull of $D/N$. If 
$R =\{\binom{d\,x}{0\,d}\mid d\in D \ \ \hbox{and} \ \ x\in E\}$ 
is the trivial extension of $D$ by $E$,  J.E. Roos proved
 that  $\lambda$-$dim (R) = n$ if and only if Krull $dim  (D) =
n$ (\cite[Theorem A']{Roo}). If $D$ is complete in its $N$-adic topology,
then $D$ is a linearly compact $D$-module, and since $E$ is
an artinian $D$-module, $R$ is a linearly compact $D$-module.
We deduce that $R$ is a local linearly compact ring, and since $R$
\ is an essential extension of a simple $R$-module, from \cite[Theorem
  7]{Anh} it follows that $R$ is a self injective ring. In the
general case,  we can prove that $R$ is self $fp$~-~injective.}
\end{remark}

\begin{corollary}
Let $R$ be an arithmetic ring of Krull dimension $0$. 

Then $\lambda$-$dim (R) \leq 2$.
\end{corollary}

\textbf{Proof.}
For every $M\in MaxSpec(R)$, any element of $MR_M$ is a zero
divisor. From \cite[Theorem 2.8]{Cou}, we deduce that $R$ is locally
self $fp$~-~injective, and from \cite[Proposition 1.2]{Cou} or
\cite[Corollary 8]{FaFa} that $R$ is self $fp$~-~injective. The result
is an immediate consequence of Theorem~\ref{T:dinj}.
\qed

\bigskip
Now, we give an example of a noncoherent Kaplansky ring $R$ with Krull
dimension $0$, which is locally coherent.

\begin{example}
\textnormal{Let $S$ be the ring defined in Proposition~\ref{P:E}.  We
  suppose that $D$ is a valuation domain with Krull dimension one and
  its maximal ideal $N$ is not finitely generated. We take $I =
  \mathbb{N}^*$.  Let $a\in N \setminus 0$ and
  $(b_n)_{n\in\mathbb{N}}$ a sequence of nonzero elements of $N$ such
  that $b_{n+1} \notin Db_n$ for every $n\in\mathbb{N}$.  We
  consider the ideal $A = Dab_0\mathbf{1} +\sum_{n\in\mathbb{N}^*}D
  ab_n\mathbf{e}_n$ and the ring $R = S/A$.  Then $(0 :_R
  \overline{a\mathbf{1}}) = R\overline{b_0\mathbf{1}} + \sum_{n\in\mathbb{N}^*} R\overline{b_n
  \mathbf{e}_n}$. Consequently $R$ is a noncoherent
ring with Krull dimension $0$.  But, for every $n\in\mathbb{N}$,
$R_{\overline{M}_{n}} \simeq D/ab_{n}D$.  Thus $R$ is locally coherent.}
\end{example}

\begin{remark} 
\textnormal{Let $S$ be the ring defined in Proposition~\ref{P:E}.  We
  suppose that $D= \mathbb{Z}_p$ (where $p$ is a prime integer) the
  ring of $p$-adic numbers and $I = \mathbb{N}^*$.  We consider $A =
\bigoplus_{n\in\mathbb{N}^*}Dp^n\mathbf{e}_n$ and $R = S/A$.  Then $R$ is
isomorphic to the example of \cite[p. 344]{Cou}. This ring is a
  Kaplansky ring which is self $fp$~-~injective, but not locally self 
$fp$~-~injective.} 
\end{remark}

\bigskip
The following proposition will be used to compute the
$\lambda$-dimension of semi-local arithmetic rings and fractionnally self
injective rings.

\begin{proposition}
\label{P:local}
Let $R$ be a valuation ring, $M$ its maximal ideal, $Z$ the subset of
zero divisors of $R$. Then the following statements are true.
\begin{enumerate}
\item $M$ is a flat module if and only if $(0 : r)$ \ is not
finitely generated for every $r\in Z \setminus 0$.

\item Let $r$ and $s$ in $R$ such that $rs \not=
0$.  Then:

\begin{itemize}
\item[i)] $(0 : rs) = ((0:r) : s)$ and $(0:r) = s(0:rs)$.
\item[ii)] If $(0 : r) \not= 0$, $(0:r)$ \ is finitely generated if
and only if $(0 : rs)$ \ is also.
\end{itemize}
\end{enumerate}
\end{proposition}

\textbf{Proof.}
\begin{enumerate}
\item Suppose that $M$ is a flat module. Let $r\in Z
\setminus 0$ and $s\in (0 : r)$.  Let $\varphi : Rr \otimes M
\rightarrow M$ be the homomorphism induced by the inclusion map $Rr
\rightarrow R$.  Then $\varphi (r \otimes s)= 0$. From
\cite[pro\-position 13 p.42]{Bou} we deduce that there exist $t\in (0 : r)$
and $m\in M$ such that $s = tm$.  Consequently $Rs \varsubsetneq (0 : r)$.

Conversely let $s\in M$ and  $r\in R$ such that $\varphi (r\otimes s)
= rs = 0$. If $r\notin Z$ then $s=0$ and $r\otimes s= 0$. If $r\in
Z\setminus 0$, since $(0 : r)$ is not finitely generated, then
there exist $t\in (0 : r)$ and $m\in M$ such that $s = tm$.  Hence 
$r\otimes s = rt\otimes m = 0$.

\item i) It is easy to get the first equality and the inclusion $s(0 :
rs) \subseteq (0:r)$.  Now, if $t\in (0 : r)$, then $t\in Rs$ 
since  $rs\not= 0$.  We deduce that there exists $c\in R$ such that
$t = cs$ and it is obvious that $c\in (0 : rs)$. Then ii) is a
consequence of i).
\qed
\end{enumerate}

\begin{theorem}
\label{T:val}
Let $R$ be a valuation ring, $M$ its
maximal ideal and $Z$ the subset of its zero divisors.  Then the
following statements are true.
\begin{enumerate}
\item If $Z = 0$ then $R$ is coherent.

\item If $Z \not= 0$ and $Z\not= M$, then 
$\lambda$-$dim (R) = 2$.

\item If $Z \not= 0$ and $Z = M$, then 
$\lambda$-$dim (R) \leq 2$.

\item If $R$ is a domain or a noncoherent ring then $M$ 
is a flat ideal and for any  $R$-module $E$ with 
$\lambda_R(E) = \infty$, $flat$-$dim (E) \leq 1$.
\end{enumerate}
\end{theorem}

\textbf{Proof.}
\begin{enumerate}
\item It is obvious.
\item We have $M = \cup_{r\in M\setminus Z}Rr$. Since 
$M$ is a direct limit of free modules, $M$ is flat. Let $E$ be
a module such that $\lambda_R (E) \geq 2$.  Then we may assume that 
$E =R/rR$, where $r\in R$.  Since
$M$ is flat, we deduce from Proposition~\ref{P:local} that $r\notin Z$ if
$r\not=  0$.  We get successively that $flat$-$dim (R/rR) \leq
1$, $\lambda_R (R/rR) \geq 3$ and  $\lambda$-$dim (R) = 2$.

\item Since $Z = M$, $R$ is self $fp$-injective by \cite[Theorem
  2.8]{Cou}. By Theorem~\ref{T:dinj}, $\lambda$-$dim (R) \leq 2$.
\item It remains to examine the case $Z = M$.  Let $r\in R$ such that 
$\lambda_R (R/rR) = \infty$.  If $R$ is not coherent, then for every
 $s\in M \setminus 0$, $(0 : s)$ is not finitely generated by
Proposition~\ref{P:local}(2). We deduce that $M$ is flat and that $r = 0$ 
or $r$ is a unit. Hence $R/rR$ is a free module.
\qed
\end{enumerate}

\begin{corollary}
\label{C:semilo}
Let $R$ be a semi-local arithmetic ring. Then
$\lambda$-$dim (R) \leq 2$ and $R$ is coherent if and only if 
$R$ is locally coherent.
\end{corollary}

\textbf{Proof.}
Since $MaxSpec(R)$ is finite, $S = \Pi_{M\in MaxSpec
(R)}R_M$ is a faithfully flat $R$-module. We deduce that\\ 
$\lambda$-$dim R \leq \lambda$-$dim S = \sup \{ \lambda$-$dim R_M \mid
M\in MaxSpec(R)\}$.
\qed

\begin{corollary}
\label{C:wdone}
Let $R$ be an arithmetic ring. We suppose
that $R_M$ is a domain or a noncoherent ring for every $M\in MaxSpec
(R)$.  Then $\lambda$-$dim (R) \leq 2$.
\end{corollary}

\textbf{Proof.}
Let $E$ be an $R$-module with $\lambda_R(E) \geq 2$.
From Theorem~\ref{T:val}, we deduce that  $flat$-$dim E_M \leq 1$ for
every $M\in MaxSpec(R)$.  Hence $flat$-$dim E \leq 1$.  We consider the
following finite 2-presentation of $E$~:

\centerline{$L_2~\overset{u_2}{\rightarrow} L_1\overset{u_1}{\rightarrow} L_0~\overset{p}{\rightarrow} E \rightarrow 0$}

Then $ker(p)$ is a finitely presented  flat $R$~-~module. We deduce
successively that $ker(p)$  , $ker(u_1)$ and $ker(u_2)$ are
finitely generated projective $R$-modules. Hence $\lambda_R (E)
\geq 3$.
\qed

\begin{corollary}
\label{C:valqu}
Let $R$ be a valuation ring and $A$ a
nonzero proper ideal of $R$. Then the following statements are true.
\begin{enumerate}
\item If $A$ is prime then $R/A$ is coherent.
\item If $A$ is finitely generated, then $R/A$ is coherent
and self $fp$-injective.
\item If $A$ is not prime and not finitely generated\\ then 
$\lambda$-$\dim (R/A) = 2$.
\end{enumerate}
\end{corollary}

\textbf{Proof.}
\begin{enumerate}
\item It is obvious.
\item We have $A = Ra$ for some $a\in R$.  If $r\notin Ra$, 
then there exists $s\notin Ra$ such that $a = rs$.  Clearly $Rs
\subseteq (Ra : r)$. Let $c\in (Ra : r)$. If $cr = 0$ then $c\in
Rs$ since $rs\not= 0$. If $cr\not= 0$, then there exists $d\in R$ such
that $cr = da =dsr$. Hence $r(c-ds) = 0$. If $ds = vc$ for some $v\in
R$, we get that $rc (1-v) = 0$. Since $rc\not= 0$, $v$ is a unit
and we obtain that $c\in Rs$, and $(Ra : r) = Rs$.

\item Since $A$ is not prime, there exist $s$ and $r\in R\setminus A$,
such that $sr\in A$.  Hence $A \varsubsetneq (A : r)$ and we prove
easily that $(A : r)$   and $(A : r)/A$ \ are not
finitely generated.
\qed
\end{enumerate}

\begin{proposition}  
\label{P:quo}
Let $R$ be an arithmetic ring and $A$ a finitely generated proper
ideal of $R$ such that $(0 : A)\subseteq J(R)$, the Jacobson radical
of $R$. Then $R/A$ \ is a coherent and self $fp$-injective ring.
\end{proposition}

\textbf{Proof.}
Then, for every maximal ideal $M$ of $R$, $AR_M$ is a
nonzero finitely generated ideal of $R_M$. By Corollary~\ref{C:valqu},
$R_M/AR_M$ is self $fp$-injective. We deduce that $R/A$ is self
$fp$-injective. 

Since in every arithmetic ring the intersection of two finitely generated
ideals is a finitely generated ideal, \cite[Corollary 1.11]{ShWi}, it is
sufficient to prove that $(A : b)$ is finitely generated for every 
$b\in R\setminus A$.  Let $M$ be a maximal ideal of $R$.  Then there
exists $a\in A$ such that $AR_M = a R_M$. Since $A$ is finitely
generated, there exists $t\in R\setminus M$ such that $AR_t =
aR_t$. Now, if $b\in aR_M$, then $b = \displaystyle{\frac{c}{s}} a$ 
for some $c\in R$ and $s\in R\setminus M$, and we get the
equality $t'sb = t'ca$ for some $t'\in R\setminus M$. We deduce
that $t's\in (A : b)$, and $(AR_{tt's} :_{R_{tt's}} b) =
t's~R_{tt's}$. If $b\notin aR_M$ then there exist $c\in R$ and
$s\notin M$ such that $a = \displaystyle{\frac{c}{s}}~b$. As in the
proof of Corollary~\ref{C:valqu}, $(R_Ma :_{R_M} b) = R_Mc$. For some
$t'\in R\setminus M$ we have $t'sa = t'cb$.  We deduce that $t'c\in (A
: b)$ and\\ $(AR_{tt's} :_{R_{tt's}} b) = t'c R_{tt's}$. Hence,
for every $M\in MaxSpec(R)$, we may assume that there exist $t_M\in
R\setminus M$ and $c_M\in (A : b)$ such that $(A R_{t_{M}}
:_{R_{t_{M}}} b) = c_MR_{t_{M}}$. A finite number of open subsets
$D(t_M)$ cover $MaxSpec(R)$.  Let $D(t_1),\ldots,D(t_n)$ be these open
subsets and $c_1,\ldots,c_n \in (A : b)$ such that 
$(A R_{t_{k}} :_{R_{t_{k}}} b) = c_k R_{t_{k}}$, for every $k$, $1\leq
k\leq n$. Then we get that $\{ c_k \mid 1\leq k\leq n\}$ generates $(A : b)$.
\qed

\begin{remark}
\textnormal{Let $S$ be the ring defined in Proposition~\ref{P:E}. We
assume that $D$ has a nonzero and nonmaximal prime ideal $J$. Let 
$B = \sum_{i\in I}J\mathbf{e}_i$, $R' = R/B$, $a\in N\setminus J$,
$b\in J\setminus 0$, $a' = a\mathbf{1}+B$ and $b' = b\mathbf{1}+B$. Then, if $R$ 
is the ring of the example~\ref{E:dim3}, we have $R =
R'/a'b'R'$. Since $\lambda$-$dim (R) = 3$, $R$ is not coherent and
not self $fp$-injective. Consequently the assumption $(0 : A)
\subseteq J(R)$ 
cannot be omitted in the Proposition~\ref{P:quo}.}
\end{remark}

\bigskip
Following V\' amos \cite{Vam},  we say that $R$ is a {\it torch ring} if
the following conditions are satisfied~:
\begin{enumerate}
\item $R$ is an arithmetical ring with at least two maximal ideals.
\item $R$ has a unique minimal prime ideal $P$ which is a
nonzero uniserial module.
\end{enumerate}

We follow T.S. Shores and R. Wiegand \cite{ShWi}, by defining a {\it canonical form} for
an $R$-module $E$ to be a decomposition $E\simeq R/I_1\oplus
R/I_2\oplus\dots\oplus R/I_n,$ where $I_1\subseteq
I_2\subseteq\dots\subseteq I_n\not= R,$ and by calling a ring $R$ a
{\it CF-ring} if every direct sum of finitely many cyclic modules has a
canonical form.
\begin{theorem}
Let $R$ be CF-ring. Then the following statements are true.
\begin{enumerate}
\item $\lambda$-$dim (R) \leq 2$.
\item $R$ is coherent if and only if $R$ is locally coherent.
\end{enumerate}
\end{theorem}

\textbf{Proof.}
In \cite[Theorem 3.12]{ShWi} it is proved that every CF-ring is
arithmetic and a finite product of indecomposable CF-rings. If $R$ is
indecomposable then $R$ is either a domain (1), or a semi-local ring
(2), or a torch ring (3). In the case (1) $R$ is coherent, and the theorem is a
consequence of Corollary~\ref{C:semilo} in the  case (2). We may
assume that $R$ is a torch ring.
Then, there is only one maximal ideal $M$ such that $P_M \not= \{
0\}$, and we have $P^2 = 0$. For every maximal ideal $N \not=
M$, $R_N$ \ is a domain. Consequently, if $R_M$ is not coherent we
deduce from Corollary~\ref{C:wdone} that $\lambda$-$dim R = 2$. Now we
assume that $R_M$ is coherent. As in the previous proposition it is
sufficient to prove that $(0 : r)$ is finitely generated for any $r\in
R$. Then we have $(0 :_{R_M}r) = sR_M$ for some $s\in R$. Since the
canonical homomorphism $R\rightarrow R_M$ is monic, $rs = 0$.

 If $r\notin P$, then $s\in P$. For every maximal ideal $N$ of $R$,
 $N\not= M$, we have  $rR_N \not= 0$ and $sR_N =0$. Consequently $(0 : r)_N = 0 = sR_N$. We deduce
that $(0 : r) = Rs$.

 If $r\in P$, then $s\notin P$ since $P^2 = 0$.
Since $R$ satisfies the condition iii) of \cite[Theorem 3.10]{ShWi},
 $V(s)$ is a finite subset of $MaxSpec(R)$. We denote $V(s) = \{ M,
 N_k \ \mid 1\leq k\leq n\}$. Since $rR_{N_{k}} = 0$ then there
 exists $s_k\notin N_k$ such that $s_kr = 0$, for every $k$, $1\leq
k\leq n$. Let $A$ be the ideal of $R$ generated by $\{ s, s_k
\ \mid \ 1 \leq k \leq n\}$.  Then for every maximal ideal $N \not= M$
we have $A_N = R_N = (0 :_{R_N} r) = (0 : r)_N$. Hence we get that $A
= (0 : r)$.
\qed  

\bigskip
In \cite{Vam} V\'amos proved that every fractionnally self-injective
ring (FSI-ring) is a CF-ring. Consequently the following corollary holds.
\begin{corollary}
Let $R$ be a fractionnally self-injective
ring. Then the following statements are true.
\begin{enumerate}
\item $\lambda$-$dim (R) \leq 2$.
\item $R$ is coherent if and only if $R$ is locally coherent.
\end{enumerate}
\end{corollary}

\section{Fractionnally self $fp$-injective rings}
\label{S:frac}
First we give a generalization of results obtained in \cite{Vam} on
fractionnally self-injective rings.

\begin{theorem}
\label{T:frac}
Let $R$ be a fractionnally self $fp$~-~injective ring. Then the
following statements are true. 
\begin{enumerate}
\item $R$ is an arithmetic ring.
\item For every proper ideal $A$ of $R$, $MinSpec (R/A)$ 
is a compact space. Moreover if $A$ is semi-prime then  $R/A$
is semihereditary.
\end{enumerate}
\end{theorem}

\textbf{Proof.}
\begin{enumerate}
\item It is the main result of  \cite{FaFa} (Theorem 1).

\item  $MinSpec (R/A)$ is homeomorphic to  $MinSpec(R/rad A)$. We may assume
that $A$ is semi-prime. Then $gl$-$w$-$dim (R/A) \leq 1$.  By
Theorem~\ref{T:redu}, $Q(R/A)$ is a Von Neumann regular ring. We deduce that 
$R/A$ is semi-hereditary from \cite[Theorem 5]{Endo} and  that $MinSpec
(R/A)$ is compact from \cite[Proposition 10]{Que}.
\qed
\end{enumerate}

\begin{remark}
\textnormal{If $R$ is the ring of our example~\ref{E:redu}, then it is
  isomorphic to the ring of \cite[Proposition 4]{FaFa} which is not
  fractionnally self $fp$~-~injective.}
\end{remark}

\bigskip
 Now we give a positive answer to a question proposed by M.
Henriksen, \cite[p. 1382]{ShoWie}. The following theorem is a generalization of
\cite[Theorem 2.4]{LLS} and \cite[Corollary 1.3]{ShoWie}.

\begin{theorem}
\label{T:H}
Every B\'{e}zout ring $R$ with compact minimal prime spectrum is Hermite.
\end{theorem}
\textbf{Proof}
Let $a$ and $b$ be in $R$.  We may assume that $a\not= 0$ and $b\not=
0$. Then there exist $a', b', d$, $m$ and $n$ in $R$ such that $a =
da'$, $b = db'$ and $ma+nb =d$.  We denote $c = ma'+nb'$, and $N$ the
nilradical of $R$. We have $(1-c)d = 0$. Since $d\not= 0$ it follows
that $c\notin N$.

First we suppose that $(N : c) = N$. Let $a^{\prime\prime},
b^{\prime\prime}, d'$, $m'$ and $n'$ be in $R$ such that $a' =
d'a^{\prime\prime}$, $b' = d'b^{\prime\prime}$ and $m'a'+n'b'
=d'$. Then $c\in Rd'$ and consequently $(N : d') = N$. Since 
$(1-m'a^{\prime\prime} - n'b^{\prime\prime})d' = 0$, 
$1-m'a^{\prime\prime} - n'b^{\prime\prime} \in N$.  Hence 
$m'a^{\prime\prime} + n'b^{\prime\prime}$ is a unit and the following
equalities hold : $a = a^{\prime\prime} d'd$, $b =
b^{\prime\prime}d'd$ and $Ra^{\prime\prime} + Rb^{\prime\prime} = R$.

Now we suppose that $(N : c) \ne N$. Since $MinSpec(R)$ is
compact $R' = R/N$ is a semi-hereditary ring. Denote $\overline r =
r+N$ for any $r\in R$. Then there exists an idempotent
$\overline{e}$ of $R'$ such that $(0 : \overline c) = R' (\overline{1}
-\overline{e})$. Since idempotents can be lifted modulo $N$, we may
assume that $e = e^2$. We deduce that $(N : c) = R(1-e)+N$, and $(1-e)
c\in N$. Let $P\in D(1-e)$. Thus $c\in P$ and $(1-c)\notin
P$. Consequently $D(1-e) \subseteq D(1-c)$, and since $(1-e)$ is an
idempotent, $(1-e) \in R(1-c)$. But $(1-c)d = 0$ and therefore $(1-e)d
= 0$ and $ed = d$. As in the proof of \cite[Theorem 2.4]{LLS} we denote
$a_1 = a'e$, $b_1 = b'e + (1-e)$, $m_1 = me$ and $n_1 = ne +
(1-e)$. Then we get $a = a_1d$, $b = b_1d$ and $m_1a+n_1b = d$.  Let
$c_1 = m_1a_1+n_1b_1 = ce + (1-e)$ and $r\in (N : c_1)$. We get
that $r(1-e) \in N$ and $rec\in N$. Hence $re \in (N : c) =
R(1-e) + N$. Consequently $re\in N$ and $(N : c_1) = N$. From the
previous part of the proof we deduce that $R$ is Hermite.
\qed

\begin{corollary}
Let $R$ be a B\'ezout ring of Krull dimension
at most one. 

If $MinSpec(R)$ is compact then $R$ is a Kaplansky ring.
\end{corollary}

\textbf{Proof.}
Let $N$ be the nilradical of $R$. By \cite[Corollary p.213]{Sho},
$R/N$ is a Kaplansky ring. From the previous theorem and 
\cite[Theorem 3]{Hen}, the result follows.
\qed 

\begin{theorem}
\label{T:fracka}
Let $R$ be a fractionnally self $fp$-injective ring . If  $R$ is a B\'ezout
ring then $R$ is a Kaplansky ring.
\end{theorem}

\textbf{Proof.}. By Theorem~\ref{T:frac} and Theorem~\ref{T:H} $R$ is
Hermite. Hence by \cite[Theorem 6]{GiHe2}, it is sufficient to prove
that for all $a$, $b$, $c\in R$ such that $Ra+Rb+Rc=R$ there
exist $p$ and $q\in R$ such that $Rpa+R(pb+qc)=R$. We put
\[A=\begin{pmatrix}a & 0\\ b & c\end{pmatrix}.\] By using the same
terminology as in \cite{LLS} let $E$ be an $R$-module named by $A$. It is easy
to check that $E$ is an $R/Rac$-module. Let $J=rad(Rac)$. It
follows that $\overline{R}=R/J$ is semihereditary by Theorem~\ref{T:frac}. Thus $\overline{E}=E/JE$ is
named by \[\overline{A}=\begin{pmatrix}\overline{a} &\overline{0}\\
\overline{b} & \overline{c}\end{pmatrix}.\] Since
  $\overline{R}$ is Hermite and $\overline{a}\cdot
  \overline{c}=\overline{0}$ we show, as in the proof of \cite[the Proposition]{Sho}, that there exist two invertible matrices
  $P$ and $Q$ and a diagonal matrix $D$ with entries in $\overline{R}$ such that $P\overline{A}Q=D$. We put \[D=\begin{pmatrix}
\overline{s} & \overline{0}\\ \overline{0} & \overline{t}
\end{pmatrix}.\] By \cite[Theorem 3.1]{LLS} we may assume that
$\overline{s}$ divides $\overline{t}$. The equality
$P^{-1}DQ^{-1}=\overline{A}$ implies that
$\overline{a},\overline{b},\overline{c}\in\overline{R}\overline{s}$.
It follows that $\overline{s}$ is a unit. Hence $\overline{E}$ is a
cyclic $\overline{R}$-module. By Nakayama Lemma it follows that $E$ is
cyclic over $R/Rac$. Hence $E$ is cyclic over $R$ too. Now we do as at
the end of the proof of \cite[Theorem 3.8]{LLS} to conclude.
\qed 

\bigskip
By \cite[Theorem 6]{FaFa} every arithmetic ring of Krull dimension zero is fractionnally self
$fp$-injective. By Theorem~\ref{T:fracka} every B\'{e}zout ring of Krull
dimension zero is Kaplansky. However it is possible to prove a more
general result. 
\begin{theorem}
Let $R$ be an arithmetic ring of Krull
dimension $0$. Then $R$ is a Kaplansky ring and an adequate
ring.
\end{theorem}
\textbf{Proof}
First we prove that $R$ is Hermite. Let  $a$ and
$b$ be in $R$. We denote $U = \{ M\in Spec(R)\mid aR_M\subseteq b R_M\}$ and
$F = \{ M\in Spec(R)\mid aR_M\not\subseteq b R_M\}$. Recall that $Spec(R)$ is
a (totally disconnected) Haussdorf compact space,  where $D(A)$ is
open and closed, for every finitely generated ideal $A$ of $R$.
Let $M\in U$.  Then there exist $c\in R$ and $t\in R\setminus M$
such that \(\displaystyle{\frac{a}{1}}= \displaystyle{\frac{cb}{t}}\).
We deduce that there
exists $s\in R\setminus M$ such that $s(ta-cb)=0$.  Hence, for
every $Q\in D(st)$, $R_Qa\subseteq R_Qb$. Consequently  $U$ is open and $F$
is closed. Now, let $M\in F$.  Then $bR_M\subset aR_M$ and there exists
$t\in R\setminus M$ such that $bR_Q\subseteq aR_Q$, for every $Q\in
D(t)$. If we denote $W_M = D(t)$, then $F \subseteq
\cup_{M\in F}W_M$. Since $F$ is compact, $F$ is
contained in a finite union $W$ of these open and closed subsets of 
$Spec(R)$.  Consequently, there exists an idempotent $e$ of $R$ 
such that 
$F \subseteq W = D(e)$ and $D(1-e) \subseteq U$. Since $R_Mb\subseteq R_Ma$
for every $M\in D(e)$, there exists $r\in R$ such that $be =
rae$.  There also exists $s\in R$ such that $a(1-e) = sb(1-e)$. Now
if we take $d = ae+b(1-e)$, $a' = s(1-e) + e$, $b' = (1-e) + re$, 
then $a = da'$, $b = db'$ and $ea' + (1-e)b' = 1$.

Now we prove that $R$ is adequate. Let $a$ and $b$ be
in $R$, $a\not= 0$.  There exists an idempotent $e$ in $R$ 
such that $D(b) = D(e)$. If we take $r = (1-e)+ae$ and $s =
a(1-e)+e$, then $a = rs$ and $Rr+Rb = R$.  Let $s'$ be a
nonunit in $R$ that divides $s$.  Then $V(s')\subseteq
V(b)$.  Hence $Rb + Rs'\not= R$.  From \cite[Theorem 8]{GiHe2}
 we deduce that $R$ is a Kaplansky ring.
\qed

\bigskip
The following proposition gives an answer to a question of \cite[p.233]{LLS}.
\begin{proposition} Let $R$ be an arithmetic ring and $J(R)$ its
Jacobson radical. Then the
following conditions are equivalent:
\begin{enumerate}
\item $R$ has a unique minimal prime ideal
\item For every $d\notin J(R)$, $(0:d)\subseteq J(R)$
\end{enumerate}
\end{proposition}
\textbf{Proof.} The implication $1\Rightarrow 2$ is easy.

$2\Rightarrow 1$. Suppose there are at least two minimal prime
ideals $I$ and $J$. Let $a\in I\setminus J$ and $P$ a maximal ideal
containing $I$. Then $IR_P$ is the nilradical of $R_P$. It follows that there exist $s\in R\setminus P$ and
a positive integer $n$ such that $sa^n=0$. Then $s\in J\setminus I$. Let $Q$
be a maximal ideal containing $J$. There also exist $t\in R\setminus
Q$ and a positive integer $m$ such that $ts^m=0$. Since $t\notin J(R)$,
$s^m\in (0:t)\subseteq J(R)$. But $s\notin P$ implies that $s^m\notin
J(R)$. Hence we get a contradiction.
\qed

\end{document}